\documentclass{article}
\usepackage{amssymb}

\usepackage{amsfonts}
\usepackage{amsmath}

\newtheorem{theorem}{Theorem}

\newtheorem{corollary}[theorem]{Corollary}

\newtheorem{definition}[theorem]{Definition}
\newtheorem{example}[theorem]{Example}

\newtheorem{lemma}[theorem]{Lemma}

\newtheorem{proposition}[theorem]{Proposition}
\newtheorem{remark}[theorem]{Remark}

\renewcommand {\Im}{\mathop{\rm Im}\nolimits}
\renewcommand {\Re}{\mathop{\rm Re}\nolimits}

\begin{document}

\title{Another look at the Burns-Krantz Theorem }
\author{David Shoikhet \\
{\small Department of Mathematics,} \ {\small ORT Braude College} {\small %
P.O. Box 78, 21982 Karmiel,}\\
{\small \ and} \ \ \\
{\small The Technion --- Israel Institute of Technology,} {\small 32000
Haifa, Israel} \\
{\small e-mails: davs27@netvision.net.il; davs@tx.technion.ac.il}}
\maketitle

\begin{abstract}
We obtain a generalization of the Burns-Krantz rigidity theorem for
holomorphic self-mappings of the unit disk in the spirit of the classical
Schwarz-Pick Lemma and its continuous version due to L.Harris via the
generation theory for one-parameter semigroups. In particular, we establish
geometric and analytic criteria for a holomorphic function on the disk with
a boundary null point to be a generator of a semigroup of linear fractional
transformations under some relations between three boundary derivatives of
the function at this point.
\end{abstract}

\bigskip

Let $\Delta $ be the open unit disk in the complex plane $%
\mathbb{C}
$ and let $Hol(\Delta ,\boldsymbol{\Omega })$ be the set of all holomorphic
functions (mappings) from $\Delta $ into $\Omega \subset
\mathbb{C}
$. In particular, the set $Hol\left( \Delta ,\Delta \right) $ of all
holomorphic self-mappings of $\Delta $ is the semigroup with respect to
composition operation.

The famous rigidity theorem of D.M.Burns and S.G.Krantz (\cite{B-K})
asserts:\medskip

$\blacklozenge $ Let $F\in Hol(\Delta ,\Delta )$ be such that

\begin{equation*}
F\left( z\right) =1+\left( z-1\right) +O\left( \left( z-1\right) ^{4}\right)
,
\end{equation*}%
as $z\rightarrow 1$. Then $F\left( z\right) \equiv z$ on $\Delta .$

Also, it was mentioned in \cite{B-K} that the exponent $4$ is sharp, and\
that it follows by the proof of the theorem that $O\left( \left( z-1\right)
^{4}\right) $ can be replaced by $o\left( \left( z-1\right) ^{3}\right) $.

R.Tauraso and F.Vlacci \cite{TR-VF} (see also, \cite{B-T-V} and \cite{TR})
established the same statement under a weaker requirement that $F\in
Hol(\Delta ,\Delta )$ satisfies the condition%
\begin{equation}
\lim_{r\rightarrow 1^{-}}\frac{F(r)-r}{(r-1)^{3}}=0.  \label{RR}
\end{equation}%
It can be shown, that, in fact, condition (\ref{RR}) is equivalent to the
condition that%
\begin{equation*}
\lim_{z\in \Delta}\Re \frac{F(z)-z}{(z-1)^{3}}=0,
\end{equation*}%
where the limit is taken in each Stolz angle (or nontangential approach
region) with vertex at the point $\ z=1$ (cf. \cite{MB}).

This result can be extended at least in three directions. In particular, we
will prove, \textit{inter alia}, the following continuous version of the
Burns-Krantz Theorem in the spirit of \cite{HL}.\medskip

$\blacklozenge $ Let $F\in Hol(\Delta ,\Delta )$ be such that%
\begin{equation*}
\mu (F):=\lim_{r\rightarrow 1^{-}}\frac{r-F(r)}{(r-1)^{3}}
\end{equation*}%
exists finitely. Then $\mu (F)$ is a nonnegative real number, and
\begin{equation*}
\left\vert F(z)-z\right\vert ^{2}\leq \mu \left( F\right) K(z),\quad z\in
\Delta ,
\end{equation*}%
for some continuous nonnegative function $K(z)$ such that
$\lim_{r\rightarrow 1^{-}}K\left(r\right) =0$. So, $F\in
Hol(\Delta ,\Delta ) $ with $F(1)=1$ and $F^{\prime }(1)$ $\left(
=\lim_{r\rightarrow 1^{+}}F^{\prime }(r)\right) =F^{\prime \prime
}(1)$ $\left( =\lim_{r\rightarrow 1^{+}}F^{\prime \prime
}(r)\right) =0$ is close to the identity mapping whenever $\mu
(F)$ is close to zero.

\bigskip To describe our approach and to present other results, we need some
notions and facts.

\bigskip $\bullet $ \ \ For a function\textbf{\ }$g\in $ $Hol$ $\left(
\Delta ,%
\mathbb{C}
\right) $ and a point $\zeta \in \partial \Delta $, we define the \emph{%
nontangential (or angular) limit }at $\zeta $ by $\angle \underset{%
z\rightarrow \zeta }{\lim }g\left( z\right) =:g\left( \zeta \right) $ if the
limit $\ \underset{z\rightarrow \zeta }{\lim }g\left( z\right) $ exists and
is the same in every \emph{nontangential approach region}%
\begin{equation*}
\Gamma (\zeta ,k):=\left\{ z\in \Delta :\left\vert z-\zeta \right\vert
<k(1-\left\vert z\right\vert )\right\} ,\quad k>1,
\end{equation*}%
with vertex at $\ \zeta \in \partial \Delta $, (see, for example, \cite%
{CCC-MBD} and \cite{SD}).

It is known (see \cite{PC-92}), that for a function $g\in Hol\left( \Delta ,%
\mathbb{C}
\right) $ and a point $\zeta \in \partial \Delta $, the angular limit $%
\angle \lim_{z\rightarrow \zeta }g^{\prime }\left( z\right) =:g^{\prime
}\left( \zeta \right) $ exists if and only if there exist the limits
\begin{equation*}
\begin{array}{ccc}
\eta =\angle \underset{z\rightarrow \zeta }{\lim }g\left( z\right) & \text{
\ \ and \ \ } & \angle \underset{z\rightarrow \zeta }{\lim }\dfrac{g\left(
z\right) -\eta }{z-\zeta }.%
\end{array}%
\end{equation*}%
In this case the last limit equals $g^{\prime }\left( \zeta \right) $.

\begin{definition}
\label{def1}For a positive integer $m$ \ we say that \ $g\in C_{A}^{m}(\zeta
),$ $\zeta \in \partial \Delta ,$\ if it \ admits the following
representation:
\end{definition}

\begin{equation}
g(z)=a_{0}+a_{1}(z-\zeta )+...+a_{m}(z-\zeta )^{m}+\gamma _{g}(z),
\label{dav2}
\end{equation}%
where

\begin{equation}
\angle \lim_{z\rightarrow \zeta }\frac{\gamma _{g}(z)}{(z-\zeta )^{m}}=0.
\label{dav3}
\end{equation}

It can be shown by induction (see, for example \cite{BP-SJ} and \cite%
{E-L-R-S}) that for $m\geq 1$ a function \ $g\in Hol(\Delta ,%
\mathbb{C}
)$ belongs to the class \ $C_{A}^{m}(\zeta ),$ if and only if there exist
the limits

\begin{equation}
\angle \lim_{z\rightarrow \zeta }g^{(k)}(z)=:g^{(k)}(\zeta )  \label{dav4}
\end{equation}%
for all \ $0\leq k\leq m$.\medskip

$\bullet $ \ The values \ $g^{(k)}(\zeta )$ \ $1\leq k\leq m$ \ are called
the \textbf{angular boundary derivatives}\emph{\ }of \ $g$ \ at the point $%
\zeta .$

So, one also writes in representation (\ref{dav2}):

\begin{equation*}
a_{k}=\frac{1}{k!}g^{(k)}(\zeta ).
\end{equation*}

\bigskip The faumous\textit{\ }Denjoy-Wolff Theorem (see, for example \cite%
{DA},\cite{WJ-26a},\cite{WJ-26b} and \cite{WJ-26c})\textit{\textbf{\ }}%
asserts:\medskip

$\blacklozenge $\ \textit{If }$F\in Hol(\Delta ,\Delta )$ \textit{is not the
identity mapping of }$\Delta $ \textit{and is not an elliptic automorphism
of }$\Delta $\textit{, then there is a unique point }$\tau \in \overline{%
\Delta }$\textit{\ that is\ an attractive point of the discrete time
semigroup }$S_{F}=\left\{ F_{n}\right\} _{n=1}^{\infty }$\textit{\ in }$%
\Delta $\textit{, defined by n-fold iterates of }$F:$ $F_{0}=I$ is the
identity mapping on $\Delta ,$ $F_{n}=F\circ F_{n-1},$ \textit{that is, }%
\begin{equation}
\underset{n\rightarrow \infty }{\lim }F_{n}\left( z\right) =\tau \text{ \ \
\ \ \ for all \ \ \ \ }z\in \Delta .  \label{D-W-point}
\end{equation}%
(see also \cite{SJH-93},\cite{CCC-MBD},\cite{RS-SD-97b}\ and \cite{SD}).

\begin{itemize}
\item This point $\tau $ is called the \textit{\textbf{Denjoy-Wolff point }}%
of\textit{\textbf{\ }}$F.$
\end{itemize}

For a point $\tau \in \overline{\Delta }$\textit{\ }$\ $\ and $%
k>1-\left\vert \tau \right\vert ^{2}$ define the sets
\begin{equation*}
D(\tau ,k)=\left\{ z\in \Delta :\frac{\left\vert 1-z\overline{\tau }%
\right\vert ^{2}}{1-\left\vert z\right\vert ^{2}}<k\right\} .
\end{equation*}

If $\tau \in \Delta $ is \textit{an interior point of} $\Delta ,$ then the
set $D(\tau ,k),$ $k>1-\left\vert \tau \right\vert ^{2},$ is exactly the
pseudo-hyperbolic disk:%
\begin{equation*}
\left\{ z\in \Delta :\frac{\left\vert z-\tau \right\vert }{\left\vert 1-z%
\overline{\tau }\right\vert }<r\right\}
\end{equation*}%
with radius $r=\sqrt{1-\frac{1-\left\vert \tau \right\vert ^{2}}{k}}.$

If $\tau \in \partial \Delta $ \textit{is a boundary point of }$\Delta ,$
then the set $D(\tau ,k),$ $k>0,$ is a horocycle internally tangent to $%
\partial \Delta $ at the point $\tau .$ Geometrically this means that $%
D(\tau ,k)$ is the disk in $\Delta $ centered at the point $\ \tau \frac{1}{%
1+k}$ with radius $\frac{k}{1+k}.$

The classical Schwarz-Pick Lemma and the boundary Schwarz-Wolff Lemma (see,
for example, \cite{SJH-93}) assert:\medskip

$\blacklozenge $ If $\tau \in \overline{\Delta }$\textit{\ is Denjoy-Wolff
point of }$F\in Hol(\Delta ,\Delta )$, \textit{then the sets }$D(\tau ,k)$%
\textit{, }$k>1-\left\vert \tau \right\vert ^{2}$,\textit{\ are }$F$-\textit{%
\ invariant, i.e.,}%
\begin{equation*}
F\left( D(\tau ,k)\right) \subset \left( D(\tau ,k)\right) .
\end{equation*}

In addition, the well-known Julia-Wolff-Carath\'{e}odory theorem
(see also \cite{SJH-93}, \cite{CCC-MBD},\ and \cite{SD})
states:\medskip

$\blacklozenge $ I\textit{f }$\tau \in \partial \Delta $\textit{\ is the
boundary Denjoy-Wolff point of }$F\in Hol(\Delta ,\Delta )$\textit{, then }$%
F $\textit{\ has at least the first order angular derivative} $F^{\prime
}\left( \tau \right) (=\angle \underset{z\rightarrow \zeta }{\lim }F^{\prime
}\left( z\right) )$ \textit{at }$\tau \in \partial \Delta ,$\textit{\ and
moreover, it is a positive real number less or equal to 1. Furthermore, if }$%
D(\tau ,k),$ $k>0,\ $\textit{is a horocycle in }$\Delta $\textit{\
internally tangent to }$\partial \Delta $\textit{\ at }$\tau \in \partial
\Delta ,$\textit{\ then }%
\begin{equation}
F\left( D(\tau ,k)\right) \subset \left( D(\tau ,\alpha k)\right) \quad
\label{inv}
\end{equation}%
\textit{for some }$\alpha \leq F^{\prime }\left( \tau \right) \leq 1.$

\textit{More generally, if }$\tau \in \partial \Delta $ \textit{is a
boundary fixed point of }$F\in Hol(\Delta ,\Delta )$\textit{, i.e., }$%
\lim_{r\rightarrow 1^{-}}F\left( r\tau \right) $\textit{\ }$=\tau $ \textit{%
(not necessarily the Denjoy-Wolff point), such that }$F^{\prime }\left( \tau
\right) (=\angle \underset{z\rightarrow \tau }{\lim }F^{\prime }\left(
z\right) )$ \textit{exists finitely, then }$F^{\prime }\left( \tau \right) $
\textit{is also a positive real number and inclusion (\ref{inv})
holds.\medskip }

$\bullet $ \ \ A fixed point $\tau \in \partial \Delta $ of $F\in Hol(\Delta
,\Delta )$ is said to be \textbf{a boundary regular fixed point }of $F$ if $%
F^{\prime }\left( \tau \right) (=\angle \underset{z\rightarrow \tau }{\lim }%
F^{\prime }\left( z\right) )$ exists finitely.

Thus, a boundary regular fixed point $\tau \in \partial \Delta $ of $F\in
Hol(\Delta ,\Delta )$ is its Denjoy-Wolff point if and only if $0<F^{\prime
}\left( \tau \right) \leq 1$. Otherwise, ($1<F^{\prime }\left( \tau \right) $%
), the point $\tau \in \partial \Delta $ is a \textbf{repelling (or
repulsive) fixed point of }$F.$

Note that, in contrast to the \textit{Denjoy-Wolff point, }a self-mapping%
\textit{\ }$F\in Hol(\Delta ,\Delta )$ may have many repelling fixed points
on $\partial \Delta .$

Following a more or less standard classification (see, for example, \cite%
{BP-SJ}) of elements in $Hol(\Delta ,\Delta )$, we observe that every
holomorphic self-mapping of $\Delta $, that is not an elliptic automorphism
of $\Delta ,$ falls in one of the following three different classes,
depending on the nature of its Denjoy-Wolff fixed point $\tau $:

\begin{itemize}
\item \textbf{Dilation type}: \ \ \ \ $\tau \in \Delta $ \ \ and \ $%
\left\vert F^{\prime }\left( \tau \right) \right\vert <1$;

\item \textbf{Hyperbolic type}: \ $\tau \in \partial \Delta $ \ and \ $%
0<F^{\prime }\left( \tau \right) <1$;

\item \textbf{Parabolic type}: \ \ \ $\tau \in \partial \Delta $ \ and \ $%
F^{\prime }\left( \tau \right) =1$.
\end{itemize}

A well-known rigidity property which follows from the Schwarz-Pick Lemma is
that \ \medskip

$\blacklozenge $ \textit{If }$\tau \in \Delta ,$\textit{\ is a fixed point
of }$F\in Hol(\Delta ,\Delta )$, \textit{then }$\left\vert F^{\prime }\left(
\tau \right) \right\vert =1$\textit{\ if and only if }$F$\textit{\ is either
the identity mapping or an elliptic automorphism of }$\Delta .$

Regarding a boundary fixed point, the following result (see \cite{TR-VF},
\cite{B-T-V} and \cite{TR}) is a modification of the Burns-Krantz rigidity
theorem \cite{B-K}.

Without loss of generality we will set in the sequel $\tau =1.$

\begin{theorem}
\label{th1}Let \ $F\in Hol(\Delta ,\Delta )\cap C_{A}^{3}(1)$ be such that \
$F(1)=F^{\prime }(1)=1$ and \ $F^{\prime \prime }(1)=F^{\prime \prime \prime
}(1)=0.$ Then $F=I,$ i.e., $F\left( z\right) =z$ for all $z\in \Delta .$
\end{theorem}

A generalization of Theorem \ref{th1} mentioned in \cite{TR-VF} (see also
\cite{C-M-P}) is the following assertion.

\begin{proposition}
\label{pr1}Let $F\in Hol(\Delta ,\Delta )$ have the Denjoy-Wolff point \ $%
\tau =1$ with \ $\alpha :=F^{\prime }(1)\left( \in (0,1]\right) .$
Assume that \ $F\in C_{A}^{3}$\bigskip $(1)$ and $\Re (F^{\prime
\prime }(1))=\alpha (\alpha -1).$ Then, the Schwarzian derivative
\begin{equation*}
S_{F}(1):=\frac{F^{\prime \prime \prime }(1)}{F^{\prime }(1)}-\frac{3}{2}%
\left( \frac{F^{\prime \prime }(1)}{F^{\prime }(1)}\right) ^{2}=0
\end{equation*}%
if and only if $F=\Phi ,$ where
\begin{equation*}
\Phi (z)=\frac{(F^{\prime \prime }(1)-2\alpha ^{2})z+\overline{F^{\prime
\prime }(1)}}{F^{\prime \prime }(1)z+\overline{F^{\prime \prime }(1)}%
-2\alpha ^{2}},\text{ \ \ \ \ }z\in \Delta ,
\end{equation*}%
is an automorphism of $\Delta $.

\textit{In particular:}

\textit{(1) }$F$\textit{\ is a hyperbolic automorphism if and only if \ }$%
\alpha \in (0,1)$\textit{\ ;}

\textit{(2) }$F$\textit{\ is a parabolic automorphism if and only if }$%
\alpha =1$\textit{\ (or equivalently,} $\Re (F^{\prime \prime }(1))=0$%
\textit{)}$\ $\textit{and }$\Im F^{\prime \prime }(1)\neq 0;$

\textit{(3)} $F$ \textit{is the identity mapping on }$\Delta $\textit{\ if
and only if }$\alpha =1$\textit{\ and }$\Im (F^{\prime \prime }(1))=0$%
\textit{\ (or equivalently, }$F^{\prime \prime }(1)=F^{\prime \prime \prime
}(1)=0$\textit{).}
\end{proposition}

\begin{remark}
\label{rem1}As a matter of fact, since $\Phi $ is a general form of
automorphisms having the Denjoy-Wolff point \ $\tau =1$ , Proposition \ref%
{pr1} can be directly reduced to\ Theorem \ref{th1} if we apply to $F$ the
automorphism $\Phi ^{-1}.$
\end{remark}

Trying to generalize Theorem \ref{th1} in another direction, one may
conjecture:\medskip

$\lozenge $\ Let \ $F\in Hol(\Delta ,\Delta )\cap C_{A}^{3}(1)$ be such that
\ $F(1)=1,$ $F^{\prime }(1)=\alpha \in (0,1]$ and \ $F^{\prime \prime
}(1)=F^{\prime \prime \prime }(1)=0.$ Then $F$ is an affine self-mapping of $%
\Delta ,$ i.e., $F\left( z\right) =\alpha z+\left( 1-\alpha \right) ,$ $z\in
\Delta .$

However, \ for $0<\alpha <1$ this conjecture is false due to the following
counterexample.

\begin{example}
\textbf{\ }\label{ex1}Let $F\in Hol(\Delta ,%
\mathbb{C}
)\cap C_{A}^{3}(1)$ be given as follows, $F\left( z\right) =\frac{1}{2}%
\left( z+1\right) +b(z-1)^{4}$. Calculations show that if $b$ is small
enough, say $b<\frac{1}{16}$, then $F$ is a self-mapping of $\Delta $ and
satisfies the hypotheses of the above conjecture with $\alpha =\frac{1}{2}.$
At the same time, $F$ is not affine..
\end{example}

$\diamondsuit $\ Nevertheless, since automorphisms of $\Delta $ as well as
affine mappings are Linear Fractional Transformations (LFT), a more general
but natural question here is concerned with\textbf{\ finding conditions on a
self-mapping }$F\in Hol(\Delta ,\Delta )\cap C_{A}^{3}(1)$\textbf{\ which
ensure that it is an LFT. \medskip }

$\diamondsuit $ \ Another question is also: \textbf{do some rigidity
properties hold if we replace the Denjoy-Wolff fixed point with any boundary
regular fixed point of }$F\in Hol(\Delta ,\Delta )$, in particular, with its
\textbf{repelling fixed point\ (}if it exists\textbf{)?\medskip }

To answer those questions, we begin with the following observation.

Let $F\in Hol(\Delta ,\Delta )$ be a holomorphic self-mapping of $\Delta $
with the boundary regular fixed point $\tau =1.$ Then, as we already
mentioned,%
\begin{equation*}
F(D(1,k))\subset D(1,\alpha k),
\end{equation*}%
where
\begin{equation*}
D(1,k)=\left\{ z\in \Delta :\frac{\left\vert 1-z\right\vert ^{2}}{%
1-\left\vert z\right\vert ^{2}}<k\right\} ,\text{ \ \ }k>0,
\end{equation*}%
is a horocycle in $\Delta $ internally tangent to $\partial \Delta $ at $%
\tau =1$ and $0<\alpha =F^{\prime }(1)$.

If \ $F$ is an LFT on $\Delta $ with the fixed point $\tau =1$, one can
result more. Namely, it can be seen that if
\begin{equation*}
a=\frac{1}{\alpha ^{2}}\left( F^{\prime \prime }(1)+\alpha (1-\alpha
)\right) ,
\end{equation*}%
then
\begin{equation*}
F(D(1,k)=D(1,\frac{\alpha k}{1+\alpha k\Re a}),
\end{equation*}%
i.e., the following equality holds%
\begin{equation}
\frac{\left\vert 1-F(z)\right\vert ^{2}}{1-\left\vert F(z)\right\vert ^{2}}=%
\frac{\alpha \left\vert 1-z\right\vert ^{2}}{(1-\left\vert
z\right\vert ^{2})+\alpha \Re a\left\vert 1-z\right\vert
^{2}},\text{ \ \ }z\in \Delta .  \label{d1}
\end{equation}%
Moreover, $F$ is an automorphism of $\Delta $ (either hyperbolic $(\alpha
\neq 1)$ or parabolic $\left( \alpha =1\right) $) if and only if $\Re %
a=0.$

It turns out that, under some smoothness conditions, equality (\ref{d1})
(and even some weaker condition ) is also sufficient for $F\in Hol(\Delta
,\Delta )$ to be linear-fractional.

\begin{theorem}
\label{th2}(cf.\cite{TR-VF}). Let $F\in Hol(\Delta ,\Delta )$ belong to $%
C_{A}^{3}(1)$ with $F\left( 1\right) (:=\lim_{r\rightarrow 1}F\left(
r\right) )=1,$ i.e.,%
\begin{equation*}
F(z)=1+F^{\prime }(1)(z-1)+\frac{1}{2}F^{\prime \prime }(1)(z-1)^{2}+\frac{1%
}{6}F^{\prime \prime \prime }(1)(z-1)^{3}+\gamma _{F}(z)
\end{equation*}%
where%
\begin{equation*}
\angle \lim_{z\rightarrow 1}\frac{\gamma _{F}(z)}{(z-1)^{3}}=0.
\end{equation*}%
Then

(A) Setting $\alpha =F^{\prime }(1)$ and $a=\frac{1}{\alpha
^{2}}(F^{\prime \prime }(1)+\alpha (1-\alpha ))$, we have $\Re
a\geq 0$;

(B) $F$ is a linear fractional transformation (LFT) on $\Delta $ if and only
if the following conditions hold:
\end{theorem}

\ \ \ \ \ \ \ \ \ \ \ \ \ \ \ \ \ \ \ \ \ \ \ \

\ \textit{(i)} $F(\Delta )\subseteq D\left( 1,\frac{1}{\Re
a}\right) $,
i.e.,%
\begin{equation*}
\frac{\left\vert 1-F(z)\right\vert ^{2}}{1-\left\vert F(z)\right\vert ^{2}}%
\leq \frac{1}{\Re a},\text{ \ \ \ \ }z\in \Delta ;
\end{equation*}%
\

\ \ \ \ \ \ \ \ \ \ \ \ \ \ \ \ \ \ \ \ \ \ \ \ \

\textit{(ii)} \ $3\left[ F^{\prime \prime }(1)\right] ^{2}=2F^{\prime \prime
\prime }(1)\cdot F^{\prime }(1)$ \textit{or , equivalently, the Schwarzian
derivative at the point} $\tau =1$%
\begin{equation*}
S_{F}(1)=\frac{F^{\prime \prime \prime }(1)}{F^{\prime }(1)}-\frac{3}{2}%
\left( \frac{F^{\prime \prime }(1)}{F^{\prime }(1)}\right) ^{2}=0.
\end{equation*}%
\textit{So, if conditions (i) and (ii) hold then for all }$z\in \Delta $%
\textit{\ we have equality (\ref{d1}).}

\textit{Moreover,\ if condition (ii) holds, then\ }$F$\textit{\ is
an automorphism of }$\Delta $\textit{\ if and only if }$\Re
a=0,$\textit{\
i.e.,}%
\begin{equation}
\Re F^{\prime \prime }(1)=\alpha (\alpha -1).  \label{AUT}
\end{equation}

\medskip \smallskip

Note that, if $\Re a=0,$ then (i) holds automatically, so Proposition %
\ref{pr1} is a direct consequence of \ Theorem \ref{th2}. Furthermore, $F$
is a parabolic automorphism if and only if $\alpha =1$. Otherwise, $\alpha
\neq 1$, $F$ is a hyperbolic automorphism; $\alpha <1,$ if and only if $\tau
=1$ is the Denjoy-Wolff (attractive) fixed point of $F;$ $\alpha >1$ if and
only if $\tau =1$ is a repelling fixed point of $F.\medskip $

\begin{remark}
\label{rem3}Note also that, if equality (\ref{AUT}) does not hold, i.e., $%
\Re a>0,$ and $\alpha >1,$ then in general $F$ must be a dilation
type self- mapping of $\Delta .$ At the same time, since $F$ is an
LFT, there is a horocycle $D$ in $\Delta ,$ $D\neq \Delta ,$ for
which $F$ is an automorhism of $D.$
\end{remark}

Indeed, if $\tau =1$ is a repelling fixed point $(\alpha >1)$,
then there is a point $\zeta \in \overline{\Delta }$, $\zeta \neq
1,$ that is the Denjoy-Wolff point of $F$. If $\zeta \in \partial
\Delta $, then $F^{\prime }(\zeta )=\frac{1}{\alpha }<1$ and $F$
is an LFT having two boundary fixed points. Hence, $F$ must be a
hyperbolic automorphism of $\Delta $ and $\Re a=0$. So, $\zeta \in
\Delta .$ It follows by (\ref{d1}) that each
horocycle $D(1,k)$, with \ $k\geq \frac{\alpha -1}{\alpha \Re a}$ is $F$%
-invariant. Moreover, $\frac{\left\vert 1-\zeta \right\vert ^{2}}{%
1-\left\vert \zeta \right\vert ^{2}}=\frac{\alpha -1}{\alpha \Re a}%
=:k_{0}$, i.e., the Denjoy-Wolff point $\zeta $ lies on the boundary of the
horocycle $D:=D(1,k_{0})$. Hence, the mapping $F$ is an automorphism of the
disk $D$.\medskip

\begin{remark}
\bigskip \label{rem4}It can be shown that, in fact, condition (i) of Theorem %
\ref{th2} implies that $S_{F}(1)$ is a real number. So condition
(ii) can be replaced by a formally weaker condition $\Re
S_{F}(1)=0.$ Moreover, \ if equality in (i) holds, then the latter
condition (hence, (ii)) is satisfied automatically (see
\cite{TR-VF}).\medskip
\end{remark}

\begin{remark}
\label{rem2}It can be shown that, by applying the Julia-Wolff-Caratheodory
Theorem applying to the function $p(z)=\frac{1+F(z)}{1-F(z)}$ (see Lemma 1
below) that condition (i) of the Theorem is equivalent to the condition:%
\begin{equation}
\frac{\left\vert 1-F(z)\right\vert ^{2}}{1-\left\vert F(z)\right\vert ^{2}}%
\leq \frac{\alpha \left\vert 1-z\right\vert ^{2}}{1-\left\vert
z\right\vert ^{2}+\alpha \Re a\left\vert 1-z\right\vert
^{2}}\text{ ,}  \label{d*}
\end{equation}%
i.e.,
\begin{equation}
F(D(1,k)\subseteq D\left( 1,\frac{\alpha k}{1+\alpha k\Re a}\right) ,%
\text{ \ \ for all \ \ }k>0.  \label{ineq}
\end{equation}
\end{remark}

In this form (where condition (i) is replaced by (\ref{d*}) or (\ref{ineq}))
Theorem \ref{th2} can also be obtained by examining the content and the
proof of Theorem 2.6 and Lemma 2.2 in \cite{TR-VF}. If $\tau =1$ is the
Denjoy-Wolff fixed point of $F,$ then this theorem can be essentially
improved. Namely, we will show that \textbf{it is enough to require that}
\textbf{condition (\ref{ineq}) holds for at least one }$k>0$\textbf{\ to
ensure that condition (ii) of Theorem \ref{th2} implies that }$F$\textbf{\
must be an LFT}. Moreover, if $F$ is of hyperbolic type mapping ($\alpha <1$%
), then one can even require a weaker inclusion than (\ref{ineq}) to obtain
the same conclusion.\medskip

\begin{theorem}
\label{TH2a}Let $F\in Hol(\Delta ,\Delta )\cap C_{A}^{3}(1)$ and let $\tau
=1 $ be the Denjoy-Wolff point of $F$ with $\alpha =F^{\prime }(1)(\in
(0,1]) $. The following assertion hold

(i) If $\lambda \in (0,1]$ and%
\begin{equation}
a_{\lambda }=\frac{1}{\alpha ^{2}}\left( \lambda F^{\prime \prime
}(1)+\alpha (1-\alpha )\right) ,  \label{1}
\end{equation}%
then $\Re a_{\lambda }\geq 0$;

(ii) $F$ is an LFT if and only if
\begin{equation}
\Re S_{F}(1):=\Re \left[ \frac{F^{\prime \prime }(1)}{F^{\prime
}(1)}-\frac{3}{2}\left( \frac{F^{\prime \prime }(1)}{F^{\prime
}(1)}\right) ^{2}\right] =0  \label{2}
\end{equation}%
and there exists $k>0$ such that%
\begin{equation}
F\left( D(1,\text{ }k\right) \subseteq D\left( 1,\text{
}\frac{k}{1+\left( k+1\right) \Re a_{\lambda }}\right) ;
\label{3}
\end{equation}%
where $a_{\lambda }$ is defined by (\ref{1}) with $\lambda =k\left(
k+1\right) ^{-1}.\medskip $
\end{theorem}

One can easily verify that
\begin{equation}
D\left( 1,\frac{\alpha k}{1+\alpha k\Re a}\right) \subseteq D\left( 1,%
\text{ }\frac{k}{1+\left( k+1\right) \Re a_{\lambda }}\right)
\label{inclusion}
\end{equation}%
whenever $0<\alpha \leq 1;$ $a_{\lambda }=\frac{1}{\alpha ^{2}}\left(
\lambda F^{\prime \prime }(1)+\alpha (1-\alpha )\right) $ and $a=\frac{1}{%
\alpha ^{2}}(F^{\prime \prime }(1)+\alpha (1-\alpha ))=a_{1\text{ }}$.
Moreover, equality in (\ref{inclusion}) holds if and only if $\alpha =1$
(i.e., $F$ is of parabolic type). Thus if $F$ is of hyperbolic type ($%
0<\alpha <1$), condition (\ref{3}) is a weaker requrement than condition (%
\ref{ineq}).\medskip

\begin{corollary}
\label{col*}Let $F\in Hol(\Delta ,\Delta )\cap C_{A}^{3}(1)$ be such that $%
F(1)=1$ and $\alpha =F^{\prime }(1)\in (0,1]$. Then $F$ is an automorphism
of $\Delta $ if and only if the following two conditions hold:

(A) $\Re S_{F}(1)=0$;

and

(B) for some $\lambda \in (0,1]$%
\begin{equation}
\lambda \Re F^{\prime \prime }(1)=\alpha (\alpha -1).  \label{Aut}
\end{equation}

Moreover, if $0<\lambda <1$ then $F$ is either the identity mapping or
parabolic automorphism of $\Delta $.
\end{corollary}

Now we point out another consequence of Theorem \ref{TH2a}, which is also a
generalization of the Burns-Krantz Theorem.\medskip

\begin{corollary}
\label{col1}Let $F\in Hol(\Delta ,\Delta )\cap C_{A}^{3}(1)$ with $0<\alpha
=F^{\prime }(1)\leq 1$. The following are equivalent:

(i) $F$ satisfies the condition%
\begin{equation*}
F(D(1,k))\subseteq D\left( 1,\frac{k\alpha }{1+k(1-\alpha )}\right) ,\text{
\ \ \ \ }k>0,
\end{equation*}%
(or, equivalently,

\begin{equation*}
F(\Delta )\subseteq D\left( 1,\frac{\alpha }{1-\alpha }\right) ,
\end{equation*}
and\
\begin{equation*}
F^{\prime \prime }(1)=F^{\prime \prime \prime }(1)=0.
\end{equation*}

(ii) $F$ is an affine mapping of the form $F(z)=\alpha z+1-\alpha .$\bigskip
\end{corollary}

\bigskip

A continuous version of the rigidity part of Theorem \ref{th2} can be given
as follows.\medskip

\begin{theorem}
\label{th3}Let \textit{\ }$F\in Hol$ $\left( \Delta ,\Delta \right) \cap
C_{A}^{3}(1)$ be such that $F(1)=1$ and $F^{\prime }(1)=\alpha (>0)$, and
let $G\in Hol(\Delta ,%
\mathbb{C}
)$ be an LFM of the form%
\begin{equation*}
G(z)=C^{-1}\left( \frac{1+z}{\alpha \left( 1-z\right) }+a\right)
\end{equation*}%
where $C(z)=\frac{1+z}{1-z}$ and $a=\frac{1}{\alpha ^{2}}(F^{\prime \prime
}(1)+\alpha (1-\alpha ))$.

The following assertions hold:

(i) $G$ is a self-mapping of $\Delta $, i.e., $\Re a\geq 0$.

(ii) if
\begin{equation*}
F(\Delta )\subset D\left( 1,\frac{1}{\Re a}\right) ,
\end{equation*}%
then the Schwarzian derivative $S_{F}(1)$ is a nonpositive real number and
the following estimate holds
\begin{equation*}
\left\vert F(z)-G(z)\right\vert ^{2}\leq -S_{F}(1)K(z)
\end{equation*}%
for some nonnegative continuous function $K(z)$.\medskip
\end{theorem}

Our approach to prove the above assertions is via rigidity properties of
generators of one-parameter continuous semigroups which are of independent
interest.\medskip

\begin{definition}
\bigskip\ \label{def2}\textit{A family }$S=\left\{ F_{t}\right\} _{t\geq
0}\subset $\textit{\ }$Hol(\Delta ,\Delta )$ \textit{is called \textbf{a}
\textbf{one-parameter continuous semigroup} on }$\Delta $\textit{\ if}
\end{definition}

\ \ \ \ \ \ \textit{(i) }$\ F_{t+s}\left( z\right) =F_{t}\left( F_{s}\left(
z\right) \right) $\textit{\ \ for all \ }$t,s\in \left[ 0,\infty \right) $%
\textit{\ and }$z\in \Delta ;$

\ \ \ \ \ \ \textit{(ii)} $\underset{t\rightarrow 0^{+}}{\lim }F_{t}\left(
z\right) =z$\textit{\ \ for all }$z\in \Delta .$

It is well known (see, for example, \cite{BE-PH} and \cite{AM-92})
that\medskip

$\blacklozenge $ \textit{The continuity condition (ii), in fact, implies,
the continuity of }$S$\textit{\ with respect to the parameter at each }$%
t\geq 0$\textit{. }Moreover, \textit{it is also differentiable in }$\left[
0,\infty \right) $\textit{\ and the limit}
\begin{equation}
\underset{t\rightarrow 0^{+}}{\lim }\frac{z-F_{t}\left( z\right) }{t}%
=:f\left( z\right) \text{, \ \ \ \ \ }z\in \Delta \text{,}
\label{generator_def}
\end{equation}%
\textit{defines a holomorphic mapping on} $\Delta $ (see also \cite{RS-SD-96}%
, \cite{RS-SD-98b} and \cite{SD}).\medskip

$\bullet $ The function $f$ \ in (\ref{generator_def}) is called the\textit{%
\textbf{\ (infinitesimal) generator of }}$S$.

\bigskip Furthermore, by using the semigroup properties, it can be shown
(see \cite{BE-PH}, \cite{RS-SD-96} and \cite{SD}) that \textit{the function }%
$u\left( t,z\right) ):=F_{t}\left( z\right) $\textit{\ is the solution of
the Cauchy problem}:%
\begin{equation}
\left\{
\begin{array}{cc}
\dfrac{\partial u\left( t,z\right) }{\partial t}+f\left( u\left( t,z\right)
\right) =0, & t\geq 0\text{,} \\
u\left( 0,z\right) =z, & z\in \Delta \text{.}%
\end{array}%
\right.  \label{cauchy}
\end{equation}

\bigskip

Denote by $\mathcal{G(}\Delta )$ the set of all holomorphic generators on $%
\Delta $. A well-known representation of $\mathcal{G}\left( \Delta \right) $
is due to E. Berkson and H. Porta \cite{BE-PH}, namely,

\begin{itemize}
\item[$\blacklozenge $] \textit{A function }$f\in Hol\left( \Delta ,%
\mathbb{C}
\right) $\textit{\ belongs to the class }$\mathcal{G}\left( \Delta \right) $%
\textit{\ if and only if there are a point }$\tau \in \overline{\Delta }$%
\textit{\ and a function }$p\in Hol\left( \Delta ,%
\mathbb{C}
\right) $ \textit{\ with }$\Re p\left( z\right) \geq 0$\textit{\
for}$\
$\textit{all} $z\in \Delta $\textit{\ such that }%
\begin{equation}
f\left( z\right) =\left( z-\tau \right) \left( 1-z\overline{\tau }\right)
p\left( z\right) \text{, \ \ \ \ }z\in \Delta \text{,}  \label{generator}
\end{equation}%
\textit{and this representation is unique.}

\textit{Moreover,}

\item[$\blacklozenge $] \textit{If }$\tau \in \Delta $\textit{\ and }$f\neq
0 $\textit{\ identically, then }$\tau $\textit{\ is the unique null point of
}$f$\textit{\ in }$\Delta $\textit{.}

\textit{If }$S=\left\{ F_{t}\right\} _{t\geq 0}\subset $ $Hol$ $\left(
\Delta ,\Delta \right) $\textit{\ is the semigroup generated by }$f$\textit{%
, then (due to the uniqueness of the solution to the Cauchy problem (\ref%
{cauchy})) }$\tau \in \Delta $ \textit{\ is a common fixed point of }$S$%
\textit{, i.e., }%
\begin{equation*}
F_{t}\left( \tau \right) =\tau \text{ \ \ \ \ \ for all \ \ \ \ }t\geq 0%
\text{.}
\end{equation*}

\item[$\blacklozenge $] \textit{In addition, if }$S$\textit{\ does not
contain an elliptic automorphism of }$\Delta $\textit{, then the point }$%
\tau \in \overline{\Delta }$\textit{\ in (\ref{generator}) is \ an
attractive point of the semigroup }$S$\textit{\ in }$\Delta $\textit{, i.e.,
}

\begin{equation}
\underset{t\rightarrow \infty }{\lim }F_{t}\left( \tau \right) =\tau \ \ \ \
\text{\ for all\ \ }\ \ z\in \Delta .  \label{D-Wp'}
\end{equation}
\end{itemize}

The last assertion is a continuous analog of the classical Denjoy-Wolff
Theorem (see \cite{AM-92},\cite{RS-SD-97b} and \cite{SD}).

$\bullet $ The point $\tau $ in (\ref{D-Wp'})\ is also called the \textit{%
\textbf{Denjoy-Wolff point of \ }}$S=\left\{ F_{t}\right\} _{t\geq 0}$%
\textbf{. \medskip }

Note that, for each $t_{0}>0$ the mappings $F_{t_{0}n},$ $n=1,2,...,$ are,
actually, iterates of the single mapping $F_{t_{0}}\in Hol(\Delta ,\Delta )$
(i.e., $F_{t_{0}n}\left( z\right) =F_{t_{0}(n-1)}\left( F_{t_{0}}(z)\right)
,\ n=1,2,...,\ F_{0}\left( z\right) =z$), so the family $S_{t_{0}}=\left\{
F_{t_{0}n}\right\} _{n=1}$ forms a discrete time semigroup with the same
Denjoy-Wolff point $\tau \in \overline{\Delta }.$

Again we observe that every semigroup $S=\left\{ F_{t}\right\} _{t\geq 0}$
of $\Delta $ generated by $f\neq 0$, that does not contain an elliptic
automorphism of $\Delta ,$ falls in one of three different classes depending
on its Denjoy-Wolff point $\tau $. By using the Berkson-Porta representation
(\ref{generator}) and an infinitesimal version of the Julia-Wolff-Carath\'{e}%
odory Theorem (see \cite{E-S1} and \cite{SD}), these classes can be
described in terms of generators as follows: if $\tau \in \overline{\Delta }$
is such that $f(\tau )=0$, then $S=\{F_{t}\}_{t\geq 0}$ is of

\begin{itemize}
\item \textbf{Dilation type}: \ \ \ \ \ $\tau \in \Delta $ \ \ and \ $\Re f^{\prime }\left( \tau \right) >0;$

\item \textbf{Hyperbolic type:} \ $\tau \in \partial \Delta $ \ and \ $%
0<f^{\prime }\left( \tau \right) <\infty ;$

\item \textbf{Parabolic type: }\ \ \ $\tau \in \partial \Delta $ \ and \ $%
f^{\prime }\left( \tau \right) =0.\medskip $
\end{itemize}

\begin{remark}
\label{rem5}It is known (see for example, \cite{RS-SD-97}) that, if $F\in
Hol(\Delta ,\Delta )$, then $f:=I-F$ belongs to $\mathcal{G}(\Delta )$; so
any criterion for $f$ to be zero implies a criterion for $F(=I-f)$ to be the
identity mapping.
\end{remark}

Another useful relation between the classes $Hol(\Delta ,\Delta )$ and $%
\mathcal{G}(\Delta )$ is the following direct consequence of the
Berkson-Porta representation formula (\ref{generator}).

\begin{itemize}
\item[$\blacklozenge $] \textit{A function }$f\in Hol\left( \Delta ,%
\mathbb{C}
\right) $\textit{\ belongs to the class }$\mathcal{G}\left( \Delta \right) $%
\textit{\ if and only if \ it admits the representation }%
\begin{equation}
f(z):=\left( z-\tau \right) \left( 1-z\overline{\tau }\right) \frac{\eta
+F\left( z\right) }{\eta -F\left( z\right) },  \label{B_P2}
\end{equation}%
\textit{with some }$\tau \in \overline{\Delta },$\textit{\ }$\eta \in
\partial \Delta $\textit{\ and }$F\in Hol(\Delta ,\Delta )$. \textit{%
Moreover, if }$\tau =\eta $\textit{, then it is a boundary regular fixed
point of }$F$\textit{\ if and only if it is the Denjoy-Wolff point for the
semigroup }$S=\left\{ F_{t}\right\} _{t\geq 0}$\textit{\ generated by }$f.$%
\textit{\ In addition, }$S$\textit{\ consists of LFT's if and only if }$F$%
\textit{\ is also an LFT.}
\end{itemize}

So, the main goal of this paper is to answer the following general
question:\medskip

$\circ $ \textbf{Given \ }$S=\left\{ F_{t}\right\} _{t\geq 0}$\textbf{\
generated by }$f\in \mathcal{G}(\Delta )$\textbf{\ with a boundary regular
common fixed point }$\tau \in \partial \Delta $\textbf{,\ find conditions on
}$f$ \textbf{which ensure that the semigroup }$S$ \textbf{consists of
linear-fractional transformations.\medskip }

\textbf{\ }Although in our last observation a mapping $F\in Hol(\Delta
,\Delta )$, (which is not the identity) having a fixed point at $\tau =1,$
may belong to any classification subclass of $Hol(\Delta ,\Delta ):$
dilation, hyperbolic or parabolic, we will see below that the function $f$
defined by (\ref{B_P2}) with $\tau =\eta (=1)$ must generate a semigroup of
hyperbolic type only. Therefore, one cannot apply this representation and
Theorem \ref{th2} to answer the question of finding rigidity conditions for
generators vanishing at the point $\tau \in \partial \Delta $, which is not
necessarily the Denjoy-Wolff point for the semigroup generated by $f$. To
answer this question we use the following notation:

\bigskip $\bullet $ A point \ $\tau \in \partial \Delta $ is said to be a
\textbf{boundary regular null point }$f$ for $f\in Hol(\Delta ,%
\mathbb{C}
)$ if the limit%
\begin{equation*}
\beta :=\angle \lim_{z\rightarrow \tau }\frac{f(z)}{z-\tau }
\end{equation*}%
exists finitely.\medskip

It was shown simultaneously in \cite{SD-01} and \cite{C-D-P}\ that, if \ $%
\tau \in \partial \Delta $ is a boundary regular null point of a generator $%
f\in \mathcal{G}(\Delta )$, then the number $\overline{\tau }\beta $ is
real. Setting again for simplicity $\tau =1$ we will complete a
characterization of the class of generators having boundary regular null
points with the following assertion.\medskip

\begin{theorem}
\label{pr2}Let \ $S=\{F_{t}\}_{t\geq 0}$ be a semigroup generated by $f\in
Hol(\Delta ,%
\mathbb{C}
)$, and let $\tau =1$. The following are equivalent: \ \ \ \ \ \ \ \ \ \ \ \
\ \ \ \ \ \ \ \ \ \ \ \ \ \ \ \ \ \ \ \ \ \ \ \ \ \ \ \ \ \ \ \ \ \ \ \ \ \
\ \ \ \ \ \ \ \ \ \ \ \ \ \ \ \ \ \ \ \ \ \ \ \ \ \ \ \ \ \ \ \ \ \ \ \ \ \
\ \ \ \ \ \ \ \ \ \ \ \ \ \ \ \ \ \ \ \ \ \ \ \ \ \ \ \ \ \ \ \ \ \ \ \ \ \
\ \ \ \ \ \ \ \ \ \ \ \ \ \ \ \ \ \ \ \ \ \ \ \ \ \ \
\end{theorem}

\textit{(i) \ }$\tau $\textit{\ is a boundary regular null point of \ }$f$%
\textit{, i.e.,}%
\begin{equation*}
\angle \lim_{z\rightarrow 1}\frac{f(z)}{z-1}=:\beta
\end{equation*}
\textit{exists finitely; \ }\ \ \ \ \ \ \ \ \ \ \ \ \ \ \ \ \ \ \ \ \ \ \ \
\ \ \ \ \ \ \ \ \ \ \ \ \ \ \ \ \ \ \ \ \ \ \ \ \ \ \ \ \ \ \ \ \ \ \ \ \ \
\ \ \ \ \ \ \ \ \ \ \ \ \ \ \ \ \ \ \ \ \ \ \ \ \ \ \ \ \ \ \ \ \ \ \ \ \ \
\ \ \ \ \ \ \ \ \ \ \ \ \ \ \ \ \ \ \ \ \ \ \ \ \ \ \ \ \ \ \ \ \ \ \ \ \ \
\ \ \ \ \ \ \ \ \ \ \ \ \ \ \ \ \ \ \ \ \ \ \ \ \ \ \ \ \ \ \ \ \

\textit{(ii) }$\tau $\textit{\ is a boundary regular common fixed point for }%
$S$\textit{, i.e., for each }$t\geq 0$\textit{\ the limit}%
\begin{equation*}
\angle \lim_{z\rightarrow 1}\frac{1-F_{t}(z)}{1-z}=:\alpha _{t}
\end{equation*}%
\textit{exists finitely. \ \ \ \ \ }\ \ \ \ \ \ \ \ \ \ \ \ \ \ \ \ \ \ \ \
\ \ \ \ \ \ \ \ \ \ \ \ \ \ \ \ \ \ \ \ \ \ \ \ \ \ \ \ \ \ \ \ \ \ \ \ \ \
\ \ \ \ \ \ \ \ \ \ \ \ \ \ \ \ \ \ \ \ \ \ \ \ \ \ \ \ \ \ \ \ \ \ \ \ \ \
\ \ \ \ \ \ \ \ \ \ \ \ \ \ \ \ \ \ \ \ \ \ \ \ \ \ \ \ \ \ \ \ \ \ \ \ \ \
\ \ \ \ \ \ \ \ \ \ \ \ \ \ \ \ \ \ \ \ \ \ \ \ \ \ \ \ \ \ \ \

(\textit{iii) there is a positive function }$\gamma :[0,\infty )\rightarrow
(0,\infty )$\textit{\ such that}%
\begin{equation*}
\frac{\left\vert 1-F_{t}(z)\right\vert ^{2}}{1-\left\vert
F_{t}(z)\right\vert ^{2}}\leq \gamma (t)\frac{\left\vert 1-z\right\vert ^{2}%
}{1-\left\vert z\right\vert ^{2}}.
\end{equation*}

\textit{(iv) }$f$\textit{\ admits the representation}%
\begin{equation*}
f(z)=-(1-z)^{2}p(z),
\end{equation*}%
\textit{where}%
\begin{equation*}
\Re p(z)\geq \frac{K}{2}\frac{1-\left\vert z\right\vert ^{2}}{%
\left\vert 1-z\right\vert ^{2}},\text{ \ \ \ }z\in \Delta ,
\end{equation*}%
\textit{for some real }$K$.

\textit{Moreover,}

\textit{(a) \ for each }$t\geq 0$\textit{\ the value }$\gamma (t)\geq
e^{-t\beta }$\textit{;}

\textit{(b) \ the maximal }$K$\textit{\ for which (iv) holds is }$\beta $%
\textit{\ , i.e.,}%
\begin{equation*}
\inf_{z\in \Delta }\frac{\left\vert 1-z\right\vert
^{2}}{1-\left\vert z\right\vert ^{2}}\cdot \Re p(z)=\frac{\beta
}{2}.
\end{equation*}%
\medskip It turns out that, in fact, condition (iv) is a criterion for \ $%
f\in Hol(\Delta ,%
\mathbb{C}
)$ with the boundary null point \ $\tau =1$ being a generator on $\Delta $.

The key for our rigidity conditions below is the following generalization of
the Berkson-Porta representation of the class $\mathcal{G}\left( \Delta
\right) .\medskip $

\begin{corollary}
\label{col2}Let \textit{\ }$f\in Hol$ $\left( \Delta ,%
\mathbb{C}
\right) $ be such that \ $\angle \lim_{z\rightarrow 1}\frac{f(z)}{z-1}=\beta
$ exists finitely.\textit{\ The function }$f$ \textit{belongs to the class }$%
\mathcal{G}\left( \Delta \right) $ if and only if it is of the form%
\begin{equation}
f(z)=-(1-z)^{2}p(z)  \label{BP2}
\end{equation}%
with
\begin{equation}
\Re p(z)\geq \frac{1}{2}K\frac{1-\left\vert z\right\vert ^{2}}{%
\left\vert 1-z\right\vert ^{2}},  \label{BP3}
\end{equation}%
for some real $K$. Moreover, $\beta $ is real and the maximal $K$ for which (%
\ref{BP3}) holds is exactly $\beta $.\medskip
\end{corollary}

\medskip Now we formulate our main regidity result.

\begin{theorem}
\label{th4}Let \textit{\ }$f\in Hol$ $\left( \Delta ,%
\mathbb{C}
\right) \cap C_{A}^{3}(1)$ with $\beta :=f^{\prime }\left( 1\right) $\
satisfy representation (\ref{BP2}) with
\begin{equation}
\inf_{z\in \Delta }\left[ \Re p(z)-\frac{1}{2}\beta
\frac{1-\left\vert z\right\vert ^{2}}{\left\vert 1-z\right\vert
^{2}}\right] =:m\geq 0. \label{NB}
\end{equation}
\end{theorem}

\textit{Then }$f$ generates a semigroup\textit{\ }$S=\{F_{t}\}_{t\geq 0}$%
\textit{\ of linear-fractional transformations if and only if the following
two conditions hold:}

\textit{(i) \ }$f^{\prime }(1)-\Re f^{\prime \prime }(1)\leq 2m$\textit{%
;}

\textit{(ii) }$f^{\prime \prime \prime }(1)=0$\textit{.}

\textit{Moreover, in this case }$m=0$\textit{\ if and only if }$f$\textit{\
is a generator of a group} \textit{of automorphisms of }$\Delta .$

\begin{remark}
\bigskip \label{rem6}We will see below that, in fact, under our setting,
\begin{equation*}
f^{\prime }(1)-\Re f^{\prime \prime }(1)=2m\geq 0.
\end{equation*}
\end{remark}

So, we have the following consequence of the above Theorem.

\begin{corollary}
\label{col3}Let $f\in \mathcal{G}(\Delta )\cap C_{A}^{3}(1)$ with $f\left(
1\right) =0$. Then $f$\textit{\ is a generator of a group }$S$ \textit{of
automorphisms of }$\Delta $ if and only if
\begin{equation}
\ f^{\prime }(1)-\Re f^{\prime \prime }(1)=\mathit{\ }f^{\prime
\prime \prime }(1)=0.\   \label{dfg}
\end{equation}%
$\ $
\end{corollary}

\textit{In addition, if (\ref{dfg}) holds, then the group }$S$\textit{\
consists of hyperbolic automorphisms if and only if }$f^{\prime }(1)\neq 0$%
\textit{. Otherwise }$(f^{\prime }(1)=0)$\textit{, \ }$S$\textit{\ consists
either of parabolic automorphisms }$(\Im f^{\prime \prime }(1)\neq 0)$%
\textit{\ \ or identity mappings (}$\Im f^{\prime \prime }(1)=0,$%
\textit{\ hence, }$f^{\prime \prime }(1)$\textit{\ equals zero).\medskip }

\begin{remark}
\label{rem7}\bigskip Using the latter assertion and Remark \ref{rem5}, we
immediately obtain Theorem \ref{th1}. Indeed, if $F\in Hol(\Delta ,\Delta
)\cap C_{A}^{3}(1)$ is such that \ $F(1)=F^{\prime }(1)=1$ and \ $F^{\prime
\prime }(1)=F^{\prime \prime \prime }(1)=0,$ then $f=I-F$ $\in \mathcal{G}%
(\Delta )\cap C_{A}^{3}(1)$ satisfies the conditions $\ f(1)=f^{\prime
}(1)=f^{\prime \prime }(1)=f^{\prime \prime \prime }(1)=0$. Hence, $f$ must
be equal identically zero, or what is one and the same, $F\left( z\right) =z$
for all $z\in \Delta $.\medskip
\end{remark}

As we have mentioned, another useful relation between self-mappings and
infinitesimal generators is representation (\ref{B_P2}). A direct
consequence of Corollary \ref{col3} and Theorem \ref{th2} (see also Lemma 6
below ) is the following assertion.\medskip

\begin{corollary}
\label{col4}Let $\ F$ be a self-mapping of $\Delta $ with the Denjoy-Wolff
point $\tau =1,$ and let $f\in \mathcal{G}(\Delta )$ be given by (\ref{B_P2}%
) with $\tau =\eta =1.$ Then $F$ is an automorphism of $\Delta $ if and only
if $\ f$ generates a group of hyperbolic automorphisms of $\Delta $.
Moreover, $F$ is parabolic if and only if \ $f^{\prime }(1)=2.\medskip $
\end{corollary}

\begin{remark}
\label{rem8}If \textit{a priory} we know that the function $p$\ in (\ref{BP2}%
) has the nonnegative real part, then the charge $\delta
_{p}(1):=\angle \lim (1-z)p(z)=\beta $. Hence, formula (\ref{BP3})
holds automatically due to Lemma 1 below, which implies that the
number $m$ defined in (\ref{NB}) equals to $\inf_{z\in \Delta }\Re
p(z)$. Consequently, in this case, our Corollary \ref{col2}
becomes the Berkson-Porta representation, of the class
$\mathcal{G}(\Delta )$ with the Denjoy-Wolff point $\tau =1$ and
condition (i) can be replaced by a formally weaker condition%
\begin{equation*}
f^{\prime }(1)-\Re f^{\prime \prime }(1)\leq 2\inf_{z\in \Delta
}\Re p(z).
\end{equation*}%
\medskip
\end{remark}

\begin{corollary}
\label{col5}Let $f\in \mathcal{G}(\Delta )\cap C_{A}^{3}(1)$ be of the form $%
f(z)=-(1-z)^{2}p(z)$. Then $f$ is a generator of a semigroup $S$ of affine
self-mappings of $\Delta $ if and only if the following two conditions hold:

(i) $\Re p(z)\geq \frac{1}{2}f^{\prime }(1)$;

(ii) $f^{\prime \prime }(1)=f^{\prime \prime \prime }(1)=0$.

In this case $\beta :=f^{\prime }(1)\geq 0$ and $f(z)=\beta (z-1)$.
\end{corollary}

\medskip The sufficient part of \ Theorem \ref{th4} can be improved as
follows.

\begin{theorem}
\label{th5}Let $f\in Hol(\Delta ,%
\mathbb{C}
)\cap C_{A}^{3}(1)$ with $\beta :=f^{\prime }(1)$ satisfy representation (%
\ref{BP2}) with (\ref{NB}), and let $g(z)=f^{\prime }(1)(z-1)+\frac{1}{2}%
f^{\prime \prime }(1)(z-1)^{2}$ be its Taylor's polynomial of degree two at $%
z=1$. The following assertions hold:

(i) $f^{\prime }(1)-\Re f^{\prime \prime }(1)\geq 0$;

(ii) $g\in \mathcal{G}(\Delta )$, i.e., $g$ generates a semigroup of LFT's
which are self-mappings of $\Delta $;

(iii) if $f^{\prime }(1)-\Re f^{\prime \prime }(1)\leq 2m,$ then $%
f^{\prime \prime \prime }(1)$ is a nonnegative real number and
\begin{equation*}
\left\vert f(z)-g(z)\right\vert ^{2}\leq \frac{1}{6}f^{\prime \prime \prime
}(1)\frac{\Re \left[ (f(z)-g(z))(1-\overline{z})^{2}\right] }{%
1-\left\vert z\right\vert ^{2}}.
\end{equation*}%
In particular, $f(z)=g(z)$ if and only if $f^{\prime \prime \prime }(1)=0$%
.\medskip
\end{theorem}

We will prove our results by using a series of simple lemmata. Since the
Berkson-Porta representation of generators contains a nonnegative real part
function, we use an approach similar to one in \cite{B-K} to analyze the
property of those functions to be linear-fractional of a certain form.
Actually, the key tool for our considerations is a modified
Julia-Wolff-Carath\'{e}odory Theorem (see Lemma \ref{lem1} bellow)
interpreted for nonnegative real part functions.

By $\Pi _{+}$ we denote the right half plane in $%
\mathbb{C}
,$ i.e., $\Pi _{+}=\left\{ z\in
\mathbb{C}
:\Re z>0\right\} .$ The well-known Riesz--Herglotz formula
\begin{equation}
p(z)=\oint\limits_{\partial \Delta }\frac{1+z\overline{\zeta }}{1-z\overline{%
\zeta }}\,dm_{p}(\zeta )+i\Im p(0)  \label{14}
\end{equation}%
establishes a linear one-to-one correspondence between class $Hol\left(
\Delta ,\overline{\Pi }_{+}\right) $ and the set of all nonnegative measure
functions $m(=m_{p})$ on the unit circle. It is easy to see that for all $%
\zeta \in \partial \Delta $ the expression $\frac{(1-z\overline{\tau })(1+z%
\overline{\zeta })}{1-z\overline{\zeta }}$ is bounded on each non-tangential
approach region at $\tau $. Then, a consequence of formula (\ref{14}) and
the Lebesgue bounded convergence theorem is that fact that for each $\tau
\in \partial \Delta $, the angular limit
\begin{equation}
\delta _{p}(\tau )=\angle \lim_{z\rightarrow \tau }(1-z\overline{\tau }%
)p(z)=2m_{p}(\tau )  \label{delta}
\end{equation}%
exists and is a nonnegative real number.

$\bullet $ The number $\delta =\delta _{p}(\tau )$ defined by (\ref{delta})
is the \textsf{charge} of the function $p\in Hol\left( \Delta ,\overline{\Pi
}_{+}\right) $ at the boundary point $\tau \in \partial \Delta $.

Denote by $C:\Delta \rightarrow \Pi _{+}$ the Cayley transform of $\Delta $
onto $\Pi _{+}:$

\begin{equation*}
C(z)=\frac{1+z}{1-z}
\end{equation*}

Applying now the Julia--Wolff-Carath\'{e}odory Theorem for the mapping $%
F:=C^{-1}\circ p\in Hol\left( \Delta ,\Delta \right) $, we get the following
assertion.\medskip

\begin{lemma}
\label{lem1}Let \ $p\in Hol(\Delta ,%
\mathbb{C}
)$ be a holomorphic function on $\Delta $ with the nonnegative real part $(%
\Re p(z)\geq 0,$ \ $z\in \Delta )$ and let $\delta _{p}(1):=\angle
\lim_{z\rightarrow 1}(1-z)p(z)$ be the charge of $p$ at $\tau =1$.
Then the
following inequality holds:%
\begin{equation}
\Re p(z)\geq \frac{1}{2}\delta _{p}(1)\frac{1-\left\vert
z\right\vert ^{2}}{\left\vert 1-z\right\vert ^{2}}.  \label{N3}
\end{equation}
\end{lemma}

\begin{lemma}
\label{lem2}(cf. \cite{M-V}) Let \ $q\in Hol(\Delta ,\overline{\Pi }_{+})$
satisfy the condition%
\begin{equation*}
\angle \lim_{z\rightarrow 1}q(z)=0.
\end{equation*}%
\ Then the angular derivative of \ $q$ at \ $z=1$ \ is either a nonpositive
real number:%
\begin{equation*}
q^{\prime }\left( 1\right) =\lim_{z\rightarrow 1}\frac{q(z)}{z-1}=:k\leq 0
\end{equation*}%
or infinity. Moreover,%
\begin{equation*}
\left\vert q(z)\right\vert ^{2}\leq -k\frac{\left\vert 1-z\right\vert ^{2}}{%
1-\left\vert z\right\vert ^{2}}\cdot \Re q(z).
\end{equation*}%
Consequently, $q^{\prime }(1)=0$ if and only if $q(z)\equiv 0$.
\end{lemma}

\textit{Proof. }If $q\neq 0$, then the function $p$ defined by $p(z):=\frac{1%
}{q(z)}$ belongs to $Hol(\Delta ,\Pi _{+})$. Noting that, $k=\frac{-1}{%
\delta _{p}(1)}$ we get our assertion from Lemma \ref{lem1}.$\blacksquare $

\begin{lemma}
\label{lem3}Let \ $p\in Hol(\Delta ,\overline{\Pi }_{+})$ admit the
following representation:%
\begin{equation}
p(z)=aC\left( z\right) +b+\gamma _{p}(z),  \label{N4}
\end{equation}%
where $\gamma _{p}(z)$ satisfies the condition%
\begin{equation}
\angle \lim_{z\rightarrow 1}\frac{\gamma _{p}(z)}{z-1}=0.  \label{N5}
\end{equation}%
Then
\end{lemma}

\textit{(i)}\ \ $a=\frac{1}{2}\delta _{p}(1)\geq 0;$\

\textit{(ii) \ the function }$\gamma _{p}(z)\equiv 0$\textit{\ if and only
if }%
\begin{equation}
\Re p\left( z\right) \geq \Re b\geq 0.  \label{N6}
\end{equation}

\textit{Proof.} Let \ $p\in Hol(\Delta ,\overline{\Pi }_{+})$ be of the form
(\ref{N4}). Then
\begin{equation*}
\delta _{p}(1)=\angle \lim_{z\rightarrow 1}(1-z)p(z)=\angle
\lim_{z\rightarrow 1}\left[ (1-z)\left( a\frac{1+z}{1-z}+b+\gamma
_{p}(z)\right) \right] =2a.
\end{equation*}%
So, assertion (i) holds. Assume that $\gamma _{p}(z)\equiv 0.$ Then $p$ is
of the form $p(z)=aC\left( z\right) +b,$ hence maps $\Delta $ into the half
plane $\left\{ z\in
\mathbb{C}
:\Re z\geq \Re b\right\} $. So, $\Re p\left( z\right) \geq
\Re b.$ At the same time, $b=p\left( z\right) -\frac{1}{2}\delta _{p}(1)%
\frac{1+z}{1-z}$; hence, by Lemma \ref{lem1}%
\begin{equation*}
\Re b=\Re p\left( z\right) -\frac{1}{2}\delta _{p}(1)\frac{%
1-\left\vert z\right\vert ^{2}}{\left\vert 1-z\right\vert ^{2}}\geq 0.
\end{equation*}%
Conversely, if (\ref{N6}) holds, then the function $p_{1}(z):=p\left(
z\right) -b$ has the nonnegative real part and has the same charge $\delta
_{p_{1}}(1)=\delta _{p}(1).$ Applying again Lemma \ref{lem1} to this
function, we obtain that $\Re \gamma _{p}(z)\geq 0.$ But by condition (%
\ref{N5}) and Lemma \ref{lem2} (see also, \cite{M-V}) this means that $%
\gamma _{p}(z)\equiv 0$.$\blacksquare $

\bigskip

If $f\in \mathcal{G}(\Delta )$, then defining $S=\{F_{t}(z)\}_{t\geq 0}$ as
the solution of the Cauchy problem (\ref{cauchy}), one can easily establish
the following fact (see, also Theorem 2.3 in \cite{B-C-D}).\medskip

\begin{lemma}
\label{lem4}Let $f\in \mathcal{G}(\Delta ).$ Then the semigroup $S\in
Hol(\Delta ,\Delta )$ generated by $f$ consists of LFT's if and only if $f$
is a polynomial of at most degree 2.\medskip
\end{lemma}

Let now \ $F\in Hol(\Delta ,\Delta )$ be a holomorphic self-mapping of $%
\Delta $ with the fixed point \ $\tau =1,$ i.e., $F(1)=1$ and \ $0<\alpha
=F^{\prime }(1)<\infty $. Then the function $f$ defined by
\begin{equation}
f(z):=-(1-z)^{2}C(F(z))  \label{N9}
\end{equation}%
belongs to $\mathcal{G}(\Delta ).$ By using the explicit form of (\ref{N9})
and direct calculations, one proves the following assertion.\medskip

\begin{lemma}
\label{lem5}The function \ $f$ defined by (\ref{N9}) belongs to class $%
C_{A}^{3}(1)$ \ if and only if \ $F$ belongs to this class. Moreover, $f$
generates a hyperbolic type semigroup with%
\begin{equation}
f^{\prime }(1)=\frac{2}{\alpha }>0,\quad \alpha =F^{\prime }(1),  \label{N10}
\end{equation}%
\begin{equation}
f^{\prime \prime }(1)=2\frac{\left[ F^{\prime }(1)\right] ^{2}-\left[
F^{\prime \prime }(1)\right] }{\left[ F^{\prime }(1)\right] ^{2}}=2\left[
\frac{1}{\alpha }-\frac{1}{\alpha ^{2}}\left[ F^{\prime \prime }(1)-\alpha
(\alpha -1)\right] \right] ,  \label{N11}
\end{equation}%
\begin{equation}
f^{\prime \prime \prime }(1)=-\frac{2}{\alpha }S_{F}(1).  \label{N12}
\end{equation}
\end{lemma}

\bigskip

We are now at the point to prove our main results.\medskip \medskip

\textbf{Proof of Theorem \ref{pr2}. }The equivalence of assertions (i) and
(ii) is proved in \cite{C-D-P} (Theorem 1) and \cite{SD-01} (Theorem 1.2).
Since by these theorems,%
\begin{equation*}
\angle \lim_{z\rightarrow 1}\frac{\partial F_{t}(z)}{\partial z}=e^{-t\beta
},\text{ \ \ \ }t\geq 0
\end{equation*}%
implication (ii)$\Rightarrow $(iii) and assertion (a) follow from the
Julia-Wolff-Carath\'{e}odory Theorem. Differentiating inequality in (iii)
with \ $\gamma (t)=e^{-t\beta }$ in $t=0^{+}$, we get
\begin{equation*}
\Re f(z)\left( \frac{\overline{z}}{1-\left\vert z\right\vert ^{2}}-%
\frac{1}{1-z}\right) =-\frac{\left\vert 1-z\right\vert
^{2}}{1-\left\vert z\right\vert ^{2}}\Re
\frac{f(z)}{(1-z)^{2}}\geq \frac{\beta }{2},
\end{equation*}%
hence, implication (iii)$\Rightarrow $(iv) follows.

Finally, assume that (iv) holds with some real $K$. Consider a holomorphic
function \ $q\in Hol(\Delta ,%
\mathbb{C}
)$ defined as follows%
\begin{equation*}
q(z)=p(z)-\frac{K}{2}\frac{1+z}{1-z}.
\end{equation*}%
Since $\Re q(z)\geq 0$, \ $z\in \Delta $, we have that the charge%
\begin{equation*}
\delta _{q}(1)=\angle \lim_{z\rightarrow 1}(1-z)q(z)\geq 0.
\end{equation*}%
In turn,%
\begin{eqnarray*}
\angle \lim_{z\rightarrow 1}\frac{f(z)}{z-1} &=&\angle \lim_{z\rightarrow
1}(1-z)p(z)=\angle \lim_{z\rightarrow 1}\left[ (1-z)q(z)+\frac{K}{2}(1+z)%
\right] = \\
&=&\delta _{q}(1)+K=:\beta .
\end{eqnarray*}%
This proves implication (iv)$\Rightarrow $(i) as well as our assertion (b).$%
\blacksquare $

\bigskip

\textbf{Proof of Corollary \ref{col2}. }If \ $f\in \mathcal{G}\left( \Delta
\right) $, then representation (\ref{BP2}) with (\ref{BP3}) follows from
Theorem \ref{pr2} and, as we mentioned above, $\beta $ must be a real
number. Conversely, assume that $f$ satisfies (\ref{BP2}) and (\ref{BP3}).
Consider a holomorphic function $g$ defined as follows:%
\begin{equation*}
g(z)=-(1-z)^{2}q(z),
\end{equation*}%
where
\begin{equation}
q\left( z\right) =p\left( z\right) -\frac{1}{2}K\frac{1+z}{1-z}.  \label{RE}
\end{equation}%
It now follows by (\ref{BP3}) that $\Re q(z)\geq 0,$ $z\in \Delta
.$ Hence, $g\in \mathcal{G}\left( \Delta \right) $ due to the
Berkson-Porta formula.

On the other hand, we have by (\ref{RE}) and (\ref{BP2}) that
\begin{equation*}
f\left( z\right) =g(z)+\frac{1}{2}K\left( z^{2}-1\right) .
\end{equation*}%
Since the second term in this formula is a generator of a one-parameter
group, we have that $\ f$ \ must belong to $\mathcal{G}\left( \Delta \right)
$.$\blacksquare $

\bigskip

\textbf{Proof of Theorem \ref{th5}.} Let \textit{\ }$f\in Hol$ $\left(
\Delta ,%
\mathbb{C}
\right) \cap C_{A}^{3}(1)$ i.e., $f$ admits the representation
\begin{equation}
f(z)=f^{\prime }(1)(z-1)+\frac{1}{2}f^{\prime \prime }(1)(z-1)^{2}++\frac{1}{%
3!}f^{\prime \prime \prime }(1)(z-1)^{3}+\gamma _{f}(z),  \notag
\end{equation}%
where%
\begin{equation}
\angle \lim_{z\rightarrow 1}\frac{\gamma _{f}(z)}{(z-1)^{3}}=0.  \label{DN3}
\end{equation}%
Then $p(z)=-f(z)(1-z)^{-2\text{ }}$ is of the form%
\begin{equation*}
p(z)=f^{\prime }(1)\frac{1}{1-z}-\frac{1}{2}f^{\prime \prime }(1)-\frac{1}{3!%
}f^{\prime \prime \prime }(1)(z-1)+\gamma _{p}(z).
\end{equation*}%
where by (\ref{DN3})
\begin{equation*}
\angle \lim_{z\rightarrow 1}\frac{\gamma _{p}(z)}{(z-1)}=-\angle
\lim_{z\rightarrow 1}\frac{\gamma _{f}(z)}{(z-1)^{3}}=0.
\end{equation*}%
We can write also:%
\begin{equation*}
p(z)=\frac{1}{2}f^{\prime }(1)\frac{1+z}{1-z}+\frac{1}{2}\left( f^{\prime
}(1)-f^{\prime \prime }(1)\right) +q(z),
\end{equation*}%
where%
\begin{equation*}
q(z)=-\frac{1}{3!}f^{\prime \prime \prime }(1)(z-1)+\gamma _{p}(z).
\end{equation*}

Noting that the function $p_{1}(z):=p(z)-\frac{\beta }{2}\frac{1+z}{1-z}$ is
of positive real part and setting
\begin{equation}
b=\frac{1}{2}\left( f^{\prime }(1)-f^{\prime \prime }(1)\right) ,  \label{N7}
\end{equation}%
we get that%
\begin{equation}
\Re b=\angle \lim_{z\rightarrow 1}\Re \left[ p_{1}(z)-q(z)\right]
=\angle \lim_{z\rightarrow 1}\Re p_{1}(z)\geq 0,  \label{N8'}
\end{equation}%
which proves assertion (i) of the theorem.

To prove assertion (ii) we note that $g(z)\left( =f^{\prime }(1)(z-1)+\frac{1%
}{2}f^{\prime \prime }(1)(z-1)^{2}\right) $ can be written in the form%
\begin{equation*}
g(z)=-b(1-z)^{2}+\frac{1}{2}\beta (z^{2}-1).
\end{equation*}%
It follows from the Berkson-Porta formula that the first term of this sum is
a generator of a parabolic type semigroup of $\Delta $, while the second
term is a generator of a group of hyperbolic automorphisms. Since $\mathcal{G%
}\left( \Delta \right) $ is a real cone (see, for example, \cite{RS-SD-98b})
we get that $g$ must belong to $\mathcal{G}\left( \Delta \right) $.
Assertion (ii) is proved.

Assume now that $\Re b\leq m$. Then $\Re q(z)\geq 0$ and
\begin{equation*}
q^{\prime }(1)=-\frac{1}{3!}f^{\prime \prime \prime }(1).
\end{equation*}%
Applying now Lemma \ref{lem2} to the function $q(z)$ and going back to the
difference $f\left( z\right) -g\left( z\right) =-\left( 1-z\right) ^{2}q(z)$%
, we complete the proof of our theorem.$\blacksquare \medskip $

\begin{remark}
\label{rem9}Note that, actually, we have from (\ref{N8'}) that $\Re %
b=\inf_{z\in \Delta }p_{1}(z)=m\geq 0$. This proves the assertion in Remark %
\ref{rem6}.
\end{remark}

Theorem \ref{th4} now is a direct consequence of Theorem \ref{th5}.\medskip

\textbf{Proof of Corollary \ref{col3}}. Denote $aut(\Delta ):=\{f\in
\mathcal{G}(\Delta ):$ such that $\ g=-f\in \mathcal{G}(\Delta )\}$. It is
clear that $f\in aut(\Delta )$ if and only if the semigroup $S$ generated by
$f$ can be extended to a one-parameter group on $\Delta $, i.e., $S$
consists of automorphisms of $\Delta $. Note that formula (\ref{dfg})
enables us to write \ $f\in \mathcal{G}(\Delta )$ in the form%
\begin{equation*}
f(z)=\frac{1}{2}\left[ f^{\prime }(1)(z^{2}-1)+i\Im f^{\prime
\prime }(1)(z-1)^{2}\right] .
\end{equation*}%
The first term of the above sum is a generator of a group of hyperbolic
automorphisms, while the second term is a generator of a group of parabolic
automorphisms. Since $aut(\Delta )$ is a real Banach algebra (see, for
example, \cite{RS-SD-98b} and \cite{A-E-R-S}), we get that\ $f\in aut(\Delta
)$. Moreover, $f$ generates a group of hyperbolic automorphisms on $\Delta $
if and only if \ $f^{\prime }(1)\neq 0.$

If \ $f^{\prime }(1)=f^{\prime \prime }(1)=f^{\prime \prime \prime }(1)=0$,
then \ $f\equiv 0$ , and following the uniqueness of the solution of the
Cauchy problem, we have that \ $F_{t}(z)\equiv z$, for all \ $t\geq 0$.

Conversely, let \ $f\in aut(\Delta ).$ It is well known (see, for example,
\cite{A-E-R-S} and \cite{SD}) that $f$ \ can be presented as
\begin{equation*}
f(z)=f(0)-\overline{f(0)}z^{2}+ikz
\end{equation*}%
with some real \ $k\in
\mathbb{R}
$ and \ $\Re f(0)\leq 0.$ Since \ $f(1)=0,$ we obtain that
\begin{equation*}
\Re f^{\prime \prime }(1)=f^{\prime }(1)=-2\Re f(0)
\end{equation*}%
and we are done.$\blacksquare $

\bigskip

\textbf{Proof} \textbf{of Theorem \ref{th2}}. Let now \ $F\in Hol(\Delta
,\Delta )\cap C_{A}^{3}(1)$ be a holomorphic self-mapping of $\Delta $ with
the fixed point \ $\tau =1,$ i.e., $F(1)=1$ and \ $0<\alpha =F^{\prime
}(1)\left( <\infty \right) .$ Then the function $f$ defined by
\begin{equation}
f(z):=-(1-z)^{2}C(F(z))
\end{equation}%
belongs to $\mathcal{G}(\Delta ).$ First we observe that \ $F$ is LFT\ if
and only if $f$ is a polynomial of at most degree 2, i.e., $f$ generates a
hyperbolic type semigroup of LFT's. By Theorem \ref{th4}, this is equivalent
that%
\begin{equation}
\Re C(F(z))\geq \frac{1}{2}(f^{\prime }(1)-\Re f^{\prime \prime
}(1))  \label{N13}
\end{equation}%
and%
\begin{equation}
f^{\prime \prime \prime }(1)=0  \label{N14}
\end{equation}%
However, by (\ref{N10}) and (\ref{N11}) we have that (\ref{N13}) can be
rewritten as
\begin{equation*}
\Re (C(F(z))\geq \frac{1}{\alpha ^{2}}\Re \left[ F^{\prime \prime
}(1)-\alpha (\alpha -1)\right] ,
\end{equation*}%
which is equivalent to inequality (\ref{d*}), hence to condition (i) of
Theorem \ref{th2}, while by (\ref{N12}) condition (\ref{N14}) is equivalent
to condition (ii) of that Theorem. Finally, note that the condition (\ref%
{AUT}) of the Theorem means that
\begin{equation*}
C(F(z))=\frac{1}{F^{\prime }(1)}\frac{1+z}{1-z}+is
\end{equation*}%
for some real $s,$ which is one and the same that $F$ is an automorphism of $%
\Delta .$\bigskip $\blacksquare $

Similarally, by using again Lemmata \ref{lem2} and \ref{lem5} and Theorem %
\ref{th5} one proves Theorem \ref{th3}.\medskip \medskip

\textbf{Proof of Theorem \ref{TH2a}. }Take any $k>0$ and consider an affine
mapping $\Phi :\Delta \rightarrow \Delta $ defined by%
\begin{equation}
\Phi (z)=\lambda z+1-\lambda ,  \label{4}
\end{equation}%
where $\lambda =\frac{k}{k+1}$. Denote $F_{1}=\Phi ^{-1}\circ F\circ \Phi $.
Since%
\begin{equation}
\Phi (\Delta )=D\left( 1,\text{ }k\right)  \label{N1'}
\end{equation}%
and%
\begin{equation}
F\left( D\left( 1,\text{ }k\right) \right) \subseteq D\left( 1,\text{ }%
\alpha k\right)  \label{N1}
\end{equation}%
we get that $F_{1}(\Delta )\subseteq \Delta $. In addition, one calculates

\begin{equation}
F_{1}^{\prime }(1)=F^{\prime }(1)=\alpha ,\text{ \ \ \ \ \ }F^{\prime \prime
}(1)=\lambda F^{\prime \prime }(1)  \label{N2'}
\end{equation}%
and%
\begin{equation}
S_{F_{1}}(1)=\lambda ^{2}S_{F}(1).  \label{5}
\end{equation}

\bigskip Furthermore, if $k>0$ is chosen such that condition (\ref{3})
holds, then we get by (\ref{N1'}) that

\begin{equation}
\Phi (F_{1}\left( \Delta \right) )=F\left( \Phi (\Delta )\right) =F\left(
D\left( 1,\text{ }k\right) \right) \subseteq D\left( 1,\text{ }\frac{k}{%
1+\left( k+1\right) \Re a_{\lambda }}\right) .  \label{Incl2}
\end{equation}%
But, for all $\mu >0$ we have
\begin{equation}
\Phi (D(1,\mu ))=D\left( 1,\frac{\mu \lambda }{1+\mu (1-\lambda )}\right)
=D\left( 1,\frac{\mu k}{\mu +k+1}\right) ,\text{ \ \ \ \ }  \label{EQ}
\end{equation}%
Setting here $\mu =\frac{1}{\Re a_{\lambda }}$ and comparing (\ref%
{Incl2}) and (\ref{EQ}) we obtain that
\begin{equation*}
F_{1}\left( \Delta \right) \subseteq D(1,\frac{1}{\Re a_{\lambda
}}).
\end{equation*}%
Now taking into account (\ref{N2'}) and (\ref{5}) and applying Theorem \ref%
{th2} we get that $F_{1},$ hence, $F(=\Phi \circ F_{1}\circ \Phi ^{-1})$ is
an LFT, and we are done.$\blacksquare \medskip $

\bigskip \textbf{Proof of Corollary \ref{col*}. }If $F$ is an automorphism
of $\Delta $ then one varifies that $S_{F}(1)=0$ and $\Re
F^{\prime \prime }(1)=F^{\prime }(1)(F^{\prime }(1)-1)$, so
condition (A) holds automatically and condition (B) holds with
$\lambda =1$.

\bigskip Converselly. If conditions (A) and (B) hold, then $F$ is an LFT by
Theorem \ref{TH2a}. Furthermore, $F_{1}=\Phi ^{-1}\circ F\circ \Phi $, where
$\Phi (z)=\lambda z+(1-\lambda )$ has the form%
\begin{equation*}
F_{1}(z)=C^{-1}\left( \frac{1}{\alpha }\cdot \frac{1+z}{1-z}+i\Im %
a_{\lambda }\right)
\end{equation*}%
where $a_{\lambda }=\lambda \Re F^{\prime \prime }(1)+\alpha
(1-\alpha ).$

Therefore, $F_{1}$ is an automorphism of $\Delta $. If $\lambda =1$, then $%
F=F_{1}$ and our assertion follows. Assume now that $\lambda <1$. Then $F$
maps the horocycle $D\left( 1,\text{ }\frac{\lambda }{1-\lambda }\right) =:D$
\ onto itself, hence $F$ is an automorphism of $D$. We clain that in this
case $\alpha =1$. Indeed, if we assume that $\alpha <1$, then $F$ must have
a fixed point $\zeta \in \partial D$, $\zeta \neq 1$. But $\partial
D\backslash \{1\}\subset \Delta $, hence $\zeta \in \Delta $ is an interion
fixed point of $\Delta $. This is a contitradiction, becouse $z=1$ is the
Denjoy-Wolff point of $F$. So, $\alpha =1$, and $F^{\prime \prime
}(1)=\alpha (\alpha -1)=0$ by (\ref{Aut}). Thus (\ref{Aut}) holds for all $%
\lambda $, in particular for $\lambda =1,$ and we are done.$\blacksquare
\medskip $

Finally, returning to the Burnz-Krants Theorem we again observe that for
each holomorphic self-mapping $F$ of $\Delta $ the function $f\in Hol(\Delta
,%
\mathbb{C}
)$ defined by $f\left( z\right) =z-F\left( z\right) $ belongs to the class $%
\mathcal{G}\left( \Delta \right) .$ Thus, applying Theorem \ref{th5} to this
function, we obtain the following assertion.\medskip

\begin{corollary}
\label{col6}Let $F\in Hol(\Delta ,\Delta )\cap C_{A}^{3}(1)$ be such that \ $%
F(1)=F^{\prime }(1)=1$ and \ $F^{\prime \prime }(1)=0.$ Then
$F^{\prime \prime \prime }(1)$ is a nonpositive real number, $\Re
\left[ (z-F(z))(1-\overline{z})^{2}\right] \geq 0$ and
\begin{equation*}
\left\vert F(z)-z\right\vert ^{2}\leq -\frac{1}{6}F^{\prime \prime
\prime }(1)\frac{\Re \left[ (z-F(z))(1-\overline{z})^{2}\right]
}{1-\left\vert z\right\vert ^{2}}
\end{equation*}%
In particular, $F(z)=z$ if and only if $F^{\prime \prime \prime }(1)=0.$
\end{corollary}

To end this discussion we conjecture another kind of rigidity
results for semigroups of holomorphic self-mappings of hyperbolic
type which is related to their asymptotic behavior at the boundary
Denjoy-Wolff point.\medskip

\textbf{Conjecture.} Let $S=\{F_{t}\}_{t\geq 0}$ be a semigroup of
hyperbolic type generated by $f\in Hol\left( \Delta ,%
\mathbb{C}
\right) $. Assume that $\tau =1$ is the Denjoy-Wolff point of $S$
and $f\in $ $\mathcal{C}_{A}^{2}(1).$ Then $f$ is a polynomial of
at most degree 2 (i.e., $S$ consists of LFT's) if and only if the
following condition holds
\begin{equation*}
\underset{t\rightarrow \infty }{\lim }\arg \left( 1-F_{t}\left( z\right)
\right) =-\arg \left( \frac{f^{\prime }(1)}{1-z}-\frac{1}{2}f^{\prime \prime
}\left( 1\right) \right) .
\end{equation*}%
Note also, that in contrast with the hyperbolic case, for semigroups of
parabolic type the latter condition is always satisfied. So, our conjecture
does not cover the parabolic case. Nevertheless, one may state the following
assertion.

$\blacklozenge $ Let $S=\{F_{t}\}_{t\geq 0}$ be a semigroup of parabolic
type generated by $f\in Hol\left( \Delta ,%
\mathbb{C}
\right) .$ Assume that $\tau =1$ is the Denjoy-Wolff point of $S$
and $f\in $ $\mathcal{C}^{3}(1).$ The semigroup
$S=\{F_{t}\}_{t\geq 0}$ consists of
parabolic automorphisms of $\Delta $ if and only if the trajectories $%
\{F_{t}(z)\}_{t\geq 0},$ $z\in \Delta ,$ converge tangentially to the point $%
\tau =1$ and $\Re F^{\prime \prime \prime }(1)=0.$

\vspace{3mm}

{\bf Acknowledgment.} The author thanks to Dr. Mark Elin for
useful discussions. Example 1 is due to him.

\bigskip

\bigskip \addcontentsline{toc}{section}{{\bf Bibliography}\hfill}

\bigskip


\markboth{\hfill Bibliography\hfill}{\hfill Bibliography\hfill}

\begin{thebibliography}{99}
\bibitem{AM-92} M. Abate, The infinitesimal generators of semigroups of
holomorphic maps, \textit{Ann. Mat. Pura Appl.} \textbf{161} (1992),
167--180.

\bibitem{A-E-R-S} D. Aharonov, M. Elin, S. {Reich,} and D. Shoikhet,
Parametric representations of semi-complete vector fields on the
unit balls in $\mathbb{C}^{n}$ and Hilbert space, \textit{Rend.
Mat. Acc. Lincei} \textbf{10} (1999), 229--253.

\bibitem{BE-PH} E. Berkson{, E.Porta, } and H. Porta, Semigroups of analytic
functions and composition operators, \textit{Michigan Math. J.} \textbf{25}
(1978), 101--115.

\bibitem{BP-SJ} P. S. Bourdon and J. H. Shapiro, Cyclic phenomena for
composition operators, \textit{Mem. Amer. Math. Soc.} \textbf{125} (1997),
no. 596.

\bibitem{B-T-V} F. Bracci, R. Tauraso and F. Vlacci, Identity principles for
commuting holomorphic self-maps of the unit disc, \textit{J. Math. Anal.
Appl.} \textbf{270} (2002), 451--473.

\bibitem{B-C-D} F. Bracci, M. D. Contreras and S. D\'{\i}az-Madrigal,
Infinitesimal generators associated with semigroups of linear fractional
maps, Journal d'Analyse Math\'{e}matique, to appear.

\bibitem{B-K} D.M. Burns and S.G. Krantz, Rigidity of holomorphic mappings
and a new Schwarz Lemma at the boundary, \textit{J. Amer. Math. Soc.}
\textbf{7} (1994), 661--676.

\bibitem{C-D-P} M. D. Contreras and S. D\'{\i}az-Madrigal, Analytic flows on
the unit disk: angular derivatives and boundary fixed points, \textit{%
Pacific J. Math.} \textbf{222} (2005), 253--286.

\bibitem{C-M-P} M. D. Contreras, S. D\'{\i}az-Madrigal and Ch. Pommerenke,
Second angular derivatives and parabolic iteration in the unit disk,
Preprint, 2006.

\bibitem{CCC-MBD} C. Cowen{\ } and B. D. MacCluer, \textit{Composition
Operators on Spaces of Analytic Functions}, CRC Press, Boca Raton, FL, 1995.

\bibitem{DA} A. Denjoy, Sur l'it\'{e}ration des fonctions analytiques,
\textit{C. R. Acad. Scie.} \textbf{182} (1926), 255--257.

\bibitem{E-L-R-S} M. Elin, M. Levinshtein, S. Reich and D. Shoikhet,
Commuting semigroups of holomorphic mappings, Preprint.

\bibitem{E-R-S2006} M. Elin, S. Reich and D. Shoikhet, A Julia--Carath\'{e}%
odory theorem for hyperbolically monotone mappings in the Hilbert ball,
Preprint.

\bibitem{E-S1} M. Elin and D. Shoikhet, Dynamic extension of the
Julia--Wolff--Carath\'{e}odory Theorem, \textit{Dynamic Systems and
Applications} \textbf{10} (2001), 421--438.

\bibitem{E-S3} M. Elin and D. Shoikhet, Semigroups with boundary fixed
points on the unit Hilbert ball and spirallike mappings, in: \textit{%
Geometric Function Theory in Several Complex Variables, 82--117, World Sci.
Publishing, River Edge, NJ}, 2004.

\bibitem{HL} L.A.Harris, A continuous form of Schwarz's Lemma in normed
linear spaces, \textit{Pacific J.Math.38} (1973), 635-639.

\bibitem{MB} T.L. Kriete and B.D. MacCluer, A rigidity theorem for
composition operators on Certain Bergman Spaces, \textit{Michigan Math. J. }%
\textbf{42} (1995), 379-386.

\bibitem{M-V} S. Migliorini and F. Vlacci, A new rigidity result for
holomorphic maps, \textit{Indag. Mathem.,N.S.} \textbf{13(4)} (2002),
537--549.

\bibitem{PC-92} Ch.~Pommerenke, \textit{Boundary Behavior of Conformal Maps}%
, Springer--Verlag, New York, Berlin, Heidelberg, 1992.

\bibitem{RS-SD-96} S. Reich and D. Shoikhet, Generation theory for
semigroups of holomorphic mappings in Banach spaces, \textit{Abstr. Appl.
Anal.} \textbf{1} (1996), 1--44.

\bibitem{RS-SD-97} S. Reich and D. Shoikhet, Semigroups and generators on
convex domains with the hyperbolic metric, \textit{Atti. Acad. Naz. Lincei}
(9) \textbf{8} (1997), 231--250.

\bibitem{RS-SD-97b} S. Reich and D. Shoikhet, The Denjoy--Wolff theorem,
\textit{Ann. Univ. Mariae Curie--Sklodowska} \textbf{51} (1997), 219--240.

\bibitem{RS-SD-98b} S. Reich and D. Shoikhet, Metric domains, holomorphic
mappings and nonlinear semigroups, \textit{Abstr. Appl. Anal.} \textbf{3}
(1998), 203--228.

\bibitem{SJH-93} J. H. Shapiro, \textit{Composition Operators and Classical
Function Theory}, Springer, Berlin, 1993.

\bibitem{SD} D. Shoikhet, \textit{Semigroups in Geometrical Function Theory}%
, Kluwer, Dordrecht, 2001.

\bibitem{SD-01} D. Shoikhet, Representations of holomorphic generators and
distortion theorems for spirallike functions with respect to a boundary
point, \textit{Int. J. Pure Appl. Math.} \textbf{5} (2003), 335--361.

\bibitem{TR} R. Tauraso, Commuting holomorphic maps of the unit disc,
\textit{Ergodic Theory Dynam. Systems} \textbf{24} (2004), 945--953.

\bibitem{TR-VF} R.Tauraso and F.Vlacci, Rigidity at the boundary for self
maps of the disk. Complex Variables:Theory and Appl., \textbf{45}%
(2),(2001),151-165.

\bibitem{WJ-26a} J. Wolff, Sur l'iteration des fonctions holomorphes dans
une region, et dont les valeurs appartiennent a cette region, \textit{C. R.
Acad. Sci.} \textbf{182} (1926), 42--43.

\bibitem{WJ-26b} J. Wolff, Sur l'iteration des fonctions bornees, \textit{C.
R. Acad. Sci.} \textbf{182} (1926), 200--201.

\bibitem{WJ-26c} J. Wolff, Sur une generalisation d'un theoreme de Schwarz,
\textit{C. R. Acad. Sci.} \textbf{182} (1926), 918--920.
\end{thebibliography}
\end{document}